\documentclass[12pt,reqno,a4wide]{amsart}

\usepackage{tikz} 
\usepackage{calrsfs}
\usepackage[charter]{mathdesign}

\allowdisplaybreaks

\newcommand{\sone}[2]{\genfrac{[}{]}{0pt}{}{#1}{#2}}

\begin{document}
\title{Logarithms of a binomial series: A Stirling number approach}

\author[Helmut Prodinger]{Helmut Prodinger}
\address{Department of Mathematics, University of Stellenbosch 7602,
	Stellenbosch, South Africa}
\email{hproding@sun.ac.za}

\begin{abstract} 
	The  $p$-th power of the logarithm of the Catalan generating function is computed using the Stirling cycle numbers.
	Instead of Stirling numbers, one may write this generating function in terms of higher order harmonic numbers.
\end{abstract}

\maketitle

\section{Introduction}

Knuth~\cite{Knuth15} proposed the exciting formula
\begin{equation*}
(\log C(z))^2=\sum_{n\ge1}\binom{2n}{n}(H_{2n-1}-H_n)\frac{z^n}{n},
\end{equation*}
where 
\begin{equation*}
C(z)=\frac{1-\sqrt{1-4z}}{2z}=\sum_{n\ge0}\frac1{n+1}\binom{2n}{n}z^n
\quad\text{and}\quad H_n=\sum_{1\le k\le n}\frac1k
\end{equation*}
with the generating function of Catalan numbers and harmonic numbers.

This formula was recently extended by Chu~\cite{Chu} to general exponents $p$. Note that Knuth talked about the exponent $1$ in his christmas lecture from 2014~\cite{christ}.

We present here a very simple approach to this question using Stirling cycle numbers; recall~\cite{GKP} that they transform falling powers into ordinary powers viz. 
\begin{equation*}
x^{\underline{n}}=\sum_{0\le k \le n}\sone{n}{k}(-1)^{n-k}x^k.
\end{equation*} 

\section{The expansion of the $p$-th power}

The substitution $z=\frac{u}{(1+u)^2}$ was presented in \cite{deBrKnRi72} and it is extremely useful when dealing with Catalan numbers and Catalan statistics. Then $C(z)=1+u$, and, by the Lagrange inversion formula~\cite{Stanley},
\begin{equation*}
u^m=\sum_{n\ge m}\frac mn\binom{2n}{n-m}z^{n}
\end{equation*}
for $m\ge1$. For $m=0$ the formula is still true when taking a limit. We now consider the bivariate generating function
\begin{equation*}
F(z,\alpha)=\sum_{p\ge0}\frac{\alpha^p}{p!}(\log C(z))^p=\exp(\alpha \log C(z))=C^\alpha(z)=(1+u)^\alpha=\sum_{m\ge0}\binom{\alpha}{m}u^m.
\end{equation*}
But
\begin{equation*}
\binom{\alpha}{m}=\frac1{m!}\alpha^{\underline{m}}=\frac1{m!}\sum_{0\le k\le m}(-1)^{m-k}\sone{m}{k}\alpha^k.
\end{equation*}
Therefore
\begin{equation*}
F(z,\alpha)=\sum_{0\le k\le m} \frac1{m!}(-1)^{m-k}\sone{m}{k}\alpha^k\sum_{n\ge m}\frac mn\binom{2n}{n-m}z^{n}.
\end{equation*}
The desired formula follows from reading off coefficients of $\alpha^p$:
\begin{equation*}
(\log C(z))^p=p![\alpha^p]F(z,\alpha)=\sum_{p \le m\le n}\frac{p!}{m!}(-1)^{m-p}\sone{m}{p}\frac mn\binom{2n}{n-m}z^{n}.
\end{equation*}

\section{Special cases}

For $p=1$, we get the instance of the Christmas lecture:
\begin{equation*}
\log C(z)=[\alpha^1]F(z,\alpha)=\sum_{1 \le m\le n}\frac{1}{m!}(-1)^{m-1}\sone{m}{1}\frac mn\binom{2n}{n-m}z^{n}.
\end{equation*}
Since $\sone{m}{1}=(m-1)!$, this leads to
\begin{equation*}
\log C(z)=[\alpha^1]F(z,\alpha)=\frac12\sum_{n\ge1}\frac1n\binom{2n}{n}z^n.
\end{equation*}

Now we turn to the instance $p=2$ from \cite{Knuth15}. (Note that $\sone m2=(m-1)!H_{m-1}$.)
\begin{align*}
2[\alpha^2]F(z,\alpha)&=\sum_{2 \le m\le n}\frac{2}{m!}(-1)^{m}\sone{m}{2}\frac mn\binom{2n}{n-m}z^{n}\\
&=2\sum_{2 \le m\le n}H_{m-1}(-1)^{m}\frac 1n\binom{2n}{n-m}z^{n}\\
&=2\sum_{1 \le j< m\le n}\frac1j(-1)^{m}\frac 1n\binom{2n}{n-m}z^{n}\\
&=2\sum_{1 \le j<  n}\frac1j(-1)^{j-1}\frac 1n\binom{2n-1}{n-j-1}z^{n}.
\end{align*}
To obtain the form proposed by Knuth, we still need to prove that
\begin{equation*}
\binom{2n}{n}(H_{2n-1}-H_n)=2\sum_{1\le j <n} \frac{(-1)^{j-1}}{j}\binom{2n-1}{n-j-1}.
\end{equation*}
Modern computer algebra systems readily simplify the difference of these two sides to 0, as expected.

\section{Connection with harmonic numbers --- the general case}

In \cite{grunberg}, there is the general formula
\begin{equation*}
	\frac1{n!}\sone{n+1}{r+1}=(-1)^r\sum_{\{r\}}
	\prod\limits_{j=1}^{l}\frac{(-1)^{i_j}}{i_j!}\bigg(\frac{H_n^{(r_j)}}{r_j}\bigg)^{i_j}.
\end{equation*}
Here, the sum is over all partitions of $r$: $r=i_1r_1+\cdots+i_lr_l$, with parts $r_1>\dots>r_l\ge1$ and positive integers
$i_1,\dots,i_l$. As an example, the partitions of $r=4$ are $4$, $3+1$, $2+ 2$, $2+1+1$, $1+1+1+1$, written alternatively as
$1\cdot4$, $1\cdot3+1\cdot1$, $2\cdot2 $, $1\cdot2+2\cdot1$, $4\cdot1$.

There   appear higher order harmonic numbers as well:
\begin{equation*}
	H_n^{(i)}=\sum_{1\le k\le n}\frac1{k^i}.
\end{equation*}
Here are the first few instances:
\begin{align*}
	\frac1{n!}\sone{n+1}{2}&=H_n,\\
	\frac1{n!}\sone{n+1}{3}&=-\frac12H_n^{(2)}+\frac12H_n^{2},\\
	\frac1{n!}\sone{n+1}{4}&=\frac13H_n^{(3)}-\frac12H_n^{(2)}H_n+\frac16H_n^3,\\
	\frac1{n!}\sone{n+1}{5}&=-\frac14H_n^{(4)}+\frac13H_n^{(3)}H_n+\frac18\big(H_n^{(2)}\big)^2-\frac14H_n^{(2)}H_n^2+\frac1{24}H_n^4.\\
\end{align*}
This allows to replace $\frac{1}{(m-1)!}\sone{m}{p}$ in
\begin{equation*}
	(\log C(z))^p=\sum_{p \le m\le n}\frac{1}{(m-1)!}\sone{m}{p}(-1)^{m-p}\frac {p!}n\binom{2n}{n-m}z^{n}
\end{equation*}
by an expression involving $H_{m-1}^{(1)},\dots,H_{m-1}^{(p-1)}$.

\section{Extension}
If instead of $u=z(1+u)^2$ we work with $u=z(1+u)^\lambda$, then we deal with the generating function
of extended (generalized) Catalan numbers
\begin{equation*}
C_\lambda(z)=\sum_{n\ge0}\binom{1+n\lambda}{n}\frac{z^n}{1+n\lambda}.
\end{equation*}

From \cite{GKP}, we infer that
\begin{equation*}
u^m=\sum_{n\ge m}\binom{\lambda n +m}{n}\frac{m}{\lambda n +m}z^n.
\end{equation*}
So
\begin{align*}
	F(z,\alpha)&=\sum_{p\ge0}\frac{\alpha^p}{p!}(\log C_\lambda(z))^p=\exp(\alpha \log C_\lambda(z))=C_\lambda^\alpha(z)=(1+u)^{\alpha}=\sum_{m\ge0}\binom{\alpha}{m}u^m\\
	&=\sum_{0\le k\le m\le n}\frac1{m!}(-1)^{m-k}\sone{m}{k}\alpha^k\binom{\lambda n +m}{n}\frac{m}{\lambda n +m}z^n.
\end{align*}
The desired formula follows from reading off coefficients of $\alpha^p$:
\begin{equation*}
	(\log C_\lambda(z))^p=p![\alpha^p]F(z,\alpha)=\sum_{p\le m\le n}
\frac{p!}{m!}(-1)^{m-p}\sone{m}{p}\binom{\lambda n +m}{n}\frac{m}{\lambda n +m}z^n.
\end{equation*}

\bibliographystyle{plain}

\begin{thebibliography}{1}
	
	
	\bibitem{deBrKnRi72}
	N.~G. De~Bruijn, D.~E. Knuth, and S.~O. Rice.
	\newblock The average height of planted plane trees.
	\newblock In R.~C. Read, editor, {\em Graph Theory and Computing}, pages
	15--22. Academic Press, 1972.
	
	
	\bibitem{GKP}
	R.~Graham, D.~E. Knuth, and O.~Patashnik.
	\newblock {\em Concrete Mathematics, second edition}.
	\newblock Addison Wesley, 1994.


	\bibitem{Knuth15}
D.~E. Knuth.
	\newblock {Log-squared of the Catalan generating function}.
	\newblock {\em Amer. Math. Monthly}, 122:390 (Problem 11832), 2015.
	
	
	\bibitem{christ}
	D.~E. Knuth.
	\newblock {$3/2$-ary trees. Annual Christmas lecture}.
	\newblock {https://www.youtube.com/watch?v=P4AaGQIo0HY},  2014.
	
	
	\bibitem{Chu}
	W.~Chu.
	\newblock Logarithms of a binomial series: Extensions of a series of Knuth.
	\newblock {\em Mathematical Communications}, 24:83--90, 2019.
	
	\bibitem{grunberg}
	D.~Gr\"{u}nberg.
	\newblock On asymptotics, {S}tirling numbers, gamma function and
	polylogs.
	\newblock {\em Results in Mathematics}, 49:89--125, 2006.


	\bibitem{Stanley}
	R.~Stanley.
	\newblock {\em Enumerative Combinatorics}.
	\newblock Wadsworth and Brooks/Cole,   1999.
	
	
\end{thebibliography}

\end{document}